\documentclass[a4paper,12pt,fleqn,leqno]{article}

\usepackage[utf8]{inputenc}
\usepackage[greek,english]{babel}
\usepackage{amssymb,amsmath,amsthm}
\usepackage{graphicx}
\usepackage[colorlinks=true,allcolors=blue]{hyperref}

\newcommand {\comment}[1] {}
\newcommand {\df}[1] {{\bfseries #1}}
\newcommand {\qt}[1] {``{#1}''}
\newcommand {\1}{{\bf1}}
\newcommand {\area} {{\textstyle\int}}

\newcommand {\C} {{\cal C}}

\newcommand {\E} {{\cal E}}

\newcommand {\enum}[2] {\{{#1},\dots,{#2}\}}

\newcommand {\floor}[1] {\lfloor#1\rfloor}
\newcommand {\id}{{\mathrm{id}}}
\newcommand {\mul} {{\:\!\circ\!\!\!\!\;\cdot\,}}
\newcommand {\N} {{\mathbb N}}

\newcommand {\Q} {{\mathbb Q}}
\newcommand {\R} {{\mathbb R}}

\newcommand {\set}[2] {\{#1\sothat#2\}}

\newcommand {\sothat} {\,:\,}
\newcommand {\sseq} {\subseteq}
\newcommand {\ti}[1] {\tilde{#1}}
\newcommand {\thet} {\vartheta}
\newcommand {\thru}[3] {{#1}_{#2},\dots,{#1}_{#3}}
\newcommand {\Z} {{\mathbb Z}}

\newcounter{substate}
\renewcommand{\thesubstate}{(\alph{substate})}
\newenvironment{substate}
{
  \begin{list}{\bf\thesubstate}
    {\usecounter{substate}
     \itemindent0em
     \settowidth\labelwidth{\bf(g)} \labelsep0.5em  
     \leftmargin\labelwidth \addtolength\leftmargin\labelsep
     \topsep0.5ex
     \itemsep0ex}
}
{\end{list}}

\newenvironment {subproof}
{
  \begin{list}{\bf\thesubstate}
    {\usecounter{substate} 
     \leftmargin0em 
     \settowidth\labelwidth{\bf(a)} \labelsep0.5em  
     \itemindent\labelwidth \addtolength\itemindent\labelsep
     \topsep0.5ex
     \itemsep0ex}
}
{\end{list}}

\newtheorem{thm}{Theorem}
\newtheorem{prp}{Proposition}

\newtheorem{obs}{Observation}
\newtheorem{cnj}{Conjecture}

\begin{document}
\title{Polynomial parametrisation of the canonical iterates
       to the solution of $-\gamma g'= g^{-1}$}
\author{Roland Miyamoto}
\maketitle
\begin{center}
\em Dedicated to the memory of
Colin Lingwood Mallows (1930--2023)\\
\
\end{center}

\begin{abstract}
\noindent
The iterates $h_0,h_1,h_2,\dotsc$
constructed in~\cite{MiSa, Mi1}
and converging to the only solution $g=h\colon[0,1]\to[0,1]$
of the iterative differential equation
$-\gamma g'= g^{-1}$, $\gamma>0$,
are parametrised by polynomials over $\Q$,
and the corresponding constant
$\gamma=\kappa\approx0.278877$
is estimated by rational numbers.
\end{abstract}

{\it MSC2020:}
  65D20, 11B83, 11Y55, 26A06, 47H10, 33E99

{\it Keywords:}
  fixed point, operator, iterative differential equation,
  special constant

\

We start by recalling the relevant definitions and facts from~\cite{MiSa,Mi1}.
For $0\leq a\leq b\leq1$
and any (Lebesgue) integrable function
$f\colon [0,1]\to\R$,
we abbreviate $\int_a^bf := \int_a^bf(x)dx$
and $\area f :=\int_0^1f$.
Let $\E$ be the set of decreasing functions
$g\colon[0,1]\to[0,1]$ satisfying $g(0)=1$ and $\area g>0$.
Given $g\in\E$,
we introduce its \df{pseudo-inverse}
$g^*\in\E$ defined by
\[
  g^*(y) := \sup g^{-1}[y,1]
  \quad\text{for } y\in[0,1],
\]
which satisfies $\area g=\area g^*$;
moreover, $g^*$ equals the compositional inverse $g^{-1}$
if $g$ is bijective.
Stretching $g\in\E$ horizontally by some factor $c=\frac1a>0$
leads to the function
$g\mul a\colon[0,c]\to[0,1]$
given by $(g\mul a)(x):=g(ax)$ for $x\in[0,c]$.

The iterative differential equation (IDE)
mentioned in the title is solved by employing
the operator $T\colon\E\to\E$ defined via
\begin{equation}\label{eqn-T}
  (Tg)(x) := \frac{\int_x^1g^*}{\area g}
  \text{ for } g\in\E \text{ and } x\in[0,1]
\end{equation}
as well as its iterations
$T^0:=\id_\E$, $T^n:=T\circ T^{n-1}$ for $n\in\N$.
We recall that $g\in\E$
solves the said IDE for some $\gamma>0$ if and only if
$g$ is a fixed point of the operator~$T$, that is, $Tg=g$,
and then $\area g=\gamma=-1/g'(0)$.
We construct a fixed point of $T$ by setting
\begin{equation}\label{eqn-hn}
  h_0:=\1_{[0,1]}\in\E \quad\text{and}\quad h_n:= T^nh_0% \in\E
  \;\;\text{for }n\in\N,
\end{equation}
and obtain as consequences that
$h_0\geq h_1\geq h_2\geq\dotsm$,
\begin{equation}\label{eqn-kappahn}
  -\kappa_nh_{n+1}' = h_n^*
  \quad \text{where} \quad
  \kappa_n := \area h_n
  \quad\text{for }n\in\N_0
\end{equation}
by~\eqref{eqn-T},
the limit
\[
  h := \lim_{n\to\infty}h_n = \inf\set{h_n}{n\in\N_0} \in\E
\]
is a fixed point of $T$,
and $h,h_1,h_2,\dotsc$ are bijective,
strictly decreasing and convex.
In~\cite{Mi1} the operator $T$ is shown to enjoy global convergence,
that is, $\lim_{n\to\infty}T^nf = h$ for every $f\in\E$.
In particular, $h$ is the only fixed point of $T$
(and of $T^n$ for every $n\in\N$),
and
\[
  \kappa := \area h = -\tfrac1{h'(0)}
  = \lim_{n\to\infty}\kappa_n
  = \inf\set{\kappa_n}{n\in\N_0}
  \approx 0.278877
\]
is the only $\gamma>0$ for which the IDE in the title
has a solution in $\E$.

\begin{figure}
\centering
\includegraphics[width=0.6\textwidth]{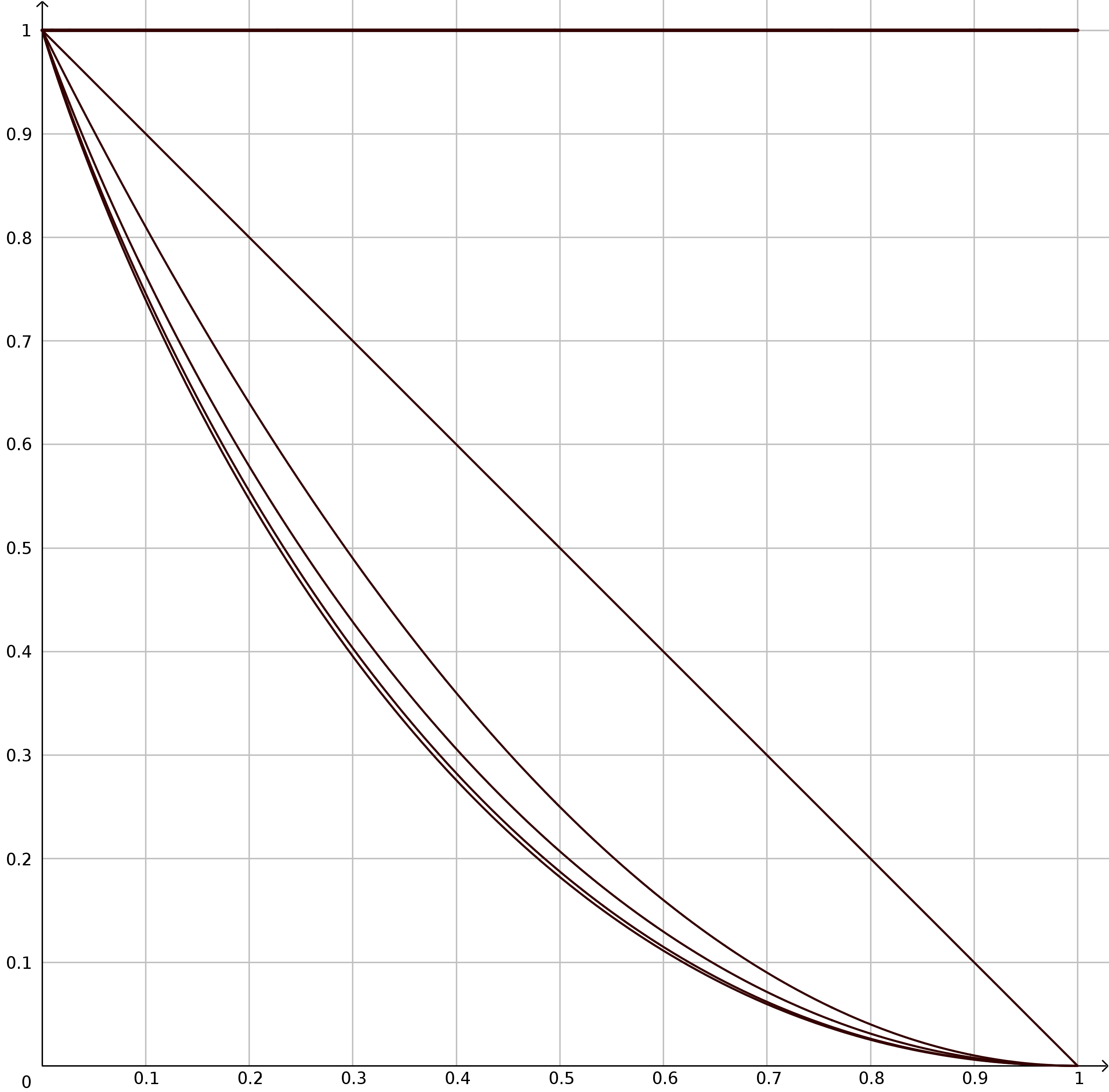}
\caption{Graphs of $\thru h05$}
\end{figure}

As observed in~\cite{Mi1},
the unique function $\ti h$ defined on the interval $[0,a]$
for some $a>0$ that satisfies the IDE
\begin{equation}\label{eqn-tih}
  -\ti h'=\ti h^*,
\end{equation}
is
$\ti h:=\tfrac1\kappa h\mul\kappa\,\colon
 [0,\tfrac1\kappa]\to[0,\tfrac1\kappa].$
Geometrically speaking,
$\ti h$ becomes its own derivative when rotated clockwise
about the origin by $90^\circ$ (into the fourth quadrant).
Owing to this property, we call $\ti h$
the \df{standard stribola} (of order~$1$;
from Greek \textgreek{στρίβω} = turn, twist)
and any function $b\cdot h\mul\frac1a\colon[0,a]\to[0,b]$
with $a,b>0$ a \df{stribola} (of order $1$),
while we refer to $h$ as the \df{unit stribola} (of order~$1$),
to $\kappa$ as the \df{stribolic constant} (of order~$1$),
to $T$ as the \df{twint} (\qt{twist and integrate})
or \df{stribolic operator} (of order~$1$),
and to the functions $h_0,h_1,h_2,\dotsc$ from~\eqref{eqn-hn}
as the \df{canonical stribolic iterates}
(of order~$1$).
(Higher order stribolas and their corresponding constants
 will be introduced and investigated in a future publication.)
In reverence to Mallows's contributions outlined below,
we also suggest the name \df{Mallows's constant} for $\kappa$.

To further motivate the calculations performed in this paper,
we want to mention some connections to combinatorics.
Colin~L.~Mallows~\cite{MaPoSl} observed
that the rows of the triangle~\cite{Le2}
underlying Levine's sequence
\[
  (l_n)_{n\in\N} = (1, 2, 2, 3, 4, 7, 14, 42, 213, 2837,
   175450, 139759600, 6837625106787,\dotsc)
\]
(see~\cite{Le1} and~\cite[p.~118]{Sl1}),
when scaled into the unit square,
seem to converge to a limit shape,
which --- its existence presumed
and after a $90^\circ$ clockwise rotation
about the point~$(\frac12,\frac12)$ ---
would necessarily equal the function $h$.
As a consequence, we would have
\[
  \lim_{n\to\infty}\frac{l_{n+1}}{l_n\cdot l_{n-1}} = \kappa,
\]
and the constants \qt{$c_1=I\approx0.277$ (last digit perhaps 6 or 8)}
conjectured by Mallows in~\cite[p.~118]{Sl1} and~\cite{Le1}
and \qt{$F(1)\dotsm\sim0.27887706$ (obtained numerically)}
conjectured by Martin Fuller in~\cite{Le2}
would also equal $\kappa$.

The OEIS entry~\cite{Ma} indicates that Mallows,
in an attempt to solve (his version of) the IDE,
must have calculated the first coefficients
of a power series for the presumed solution.
In our terminology, let $\tau=\ti h(\tau)$ be the fixpoint
of the standard stribola $\ti h$ and $b_n:=\ti h^{(n)}(\tau)$
the $n$-th derivative of $\ti h$ at $\tau$.
Then Faà di Bruno's formula applied to~\eqref{eqn-tih}
yields the identities
\[
  \sum_{k=1}^n B_{n,k}(\thru b1{n-k+1})\cdot b_{k+1} = 0
  \quad\text{for } n = 2,3,4,\dotsc
\]
where $B_{n,k}$ are the exponential Bell polynomials.
These identities can be used to express all the $b_n$
in terms of~$\tau$ as is made more explicit in~\cite{Ma},
namely
\begin{align*}
  b_0 &= \tau, \quad
  b_1 = -\tau, \quad
  b_2 = \tau^{-1}, \quad
  b_3 = -\tau^{-4}, \quad
  b_4 = 3\tau^{-7} - \tau^{-8}, \\
  b_5 &= -15\tau^{-10} + 10\tau^{-11}
        - 3\tau^{-12} + \tau^{-13},\\
  b_6 &= 105\tau^{-13} - 105\tau^{-14} + 55\tau^{-15}
      - 30\tau^{-16} + 10\tau^{-17} - 3\tau^{-18} + \tau^{-19},
  \quad\dotsc
\end{align*}

Concerning the canonical stribolic iterates
defined by~\eqref{eqn-hn}, we note that
$h_0^*=h_0$, $h_1(x)=1-x=h_1^*(x)$,
$h_2(x)=(1-x)^2$, $h_2^*(x)=1-\sqrt x$
and $h_3(x) = 1-3x+2x^{\frac32}$
for $x\in[0,1]$ and thereby
$\kappa_0=1$, $\kappa_1=\tfrac12$,
$\kappa_2=\tfrac13$ and $\kappa_3=\tfrac3{10}$.
We cannot expect explicit function terms
for $h_4,h_5,\dotsc$,
but instead are aiming at parametrising them.
To this end, we introduce the compositions
\[
  q_n := h_n\circ\dotsm\circ h_1\colon[0,1]\to[0,1]
  \quad\text{for }n\in\N_0,
\]
where, as usual, we agree that the empty composition
is $q_0=\id_{[0,1]}$.
For $n\in\N$,
the equation $y=h_n(x)$ is obviously parametrised
by $x=q_{n-1}(t)$, $y=q_n(t)$, $t\in[0,1]$.
Because $X:=\id_{[0,1]}$ is transcendental over $\R$,
we can look at the univariate polynomial ring $\R[X]$
in both ways,
algebraically and as a subspace of $\C^0[0,1]$,
the space of continuous functions on $[0,1]$.
We consider the discrete valuation
$v_X\colon\R[X]\to\N_0\cup\{\infty\}$
associated with its prime element $X$ (cf.~\cite[p.~4f]{St})
and the degree map $\deg\colon\R[X]\to\N_0\cup\{-\infty\}$.
The following theorem reveals
that the $q_n$ are polynomials over $\Q$
% the first few of which are displayed in Table~\ref{tbl-kapqn},
and their degrees are given by the Fibonacci sequence
$(F_0,F_1,F_2,\dotsc) := (0,1,1,2,3,5,8,13,\dotsc)$.

\begin{thm}\label{thm-qn}
	Let $n\in\N$ and
	$Q_n(x):=\int_0^xq_n'q_{n-1}$ for $x\in[0,1]$.
	Then
	\begin{substate}
	\item $q_n = h_n\circ q_{n-1} \colon[0,1]\to[0,1]$
	  is bijective and continuously differentiable
	  with $q_n(0)= n\bmod2\in\{0,1\}$
	  and $q_n(1)-q_n(0)=(-1)^n$,
	\item $\kappa_n = (-1)^nQ_n(1)$ and
	  $q_{n+1} = q_{n+1}(0) - \frac1{\kappa_n}Q_n$,
	\item $\kappa_n\in\Q$ and $q_n\in\Q[X]$ with
	  $v_X(q_n-q_n(0)) = 2^{\floor{\frac n2}}$
	  and $\deg q_n = F_{n+1}$;
	  moreover, there exists $d_n\in\N$ such that
	  $d_n\cdot(q_n-q_n(0))\in\Z[X]$
	  is primitive.
	\end{substate}
	Wherever it makes sense, the above statements also hold for $n=0$.
\end{thm}
\begin{proof}
  \begin{subproof}
  \item follows directly from the definition of $q_n$.
    Clearly, $q_0=\id_{[0,1]}=X$ is also
    continuously differentiable and bijective
    with $q_0(0)=0$ and $q_0(1)=1$.
  \item Equation~\eqref{eqn-kappahn} implies
    $h_{n+1}'(h_n(x)) = -\frac x{\kappa_n}$ for $x\in[0,1]$,
    hence
    \[
      q_{n+1}' = (h_{n+1}\circ q_n)' = q_n'\cdot(h_{n+1}'\circ q_n)
      = q_n'\cdot(h_{n+1}'\circ h_n\circ q_{n-1})
      = -\tfrac1{\kappa_n}q_n'q_{n-1}
    \]
    by (a), and integration yields
    $q_{n+1}(x) = q_{n+1}(0)-\frac1{\kappa_n}Q_n(x)$ for $x\in[0,1]$.
    The first assertion follows by substituting $x=1$
    and using~(a) again.
  \item By definition, $\kappa_0=1$, $q_0=X$ and $q_1=1-X$.
    Using~(b) and induction, we obtain
    $\kappa_n\in\Q$ and $q_n,Q_n\in\Q[X]$.
    From~(b) we also conclude that
    $\deg q_{n+1} = \deg Q_n = 1 + \deg(q_n'q_{n-1})
                  = \deg q_n + \deg q_{n-1}$,
    thus $\deg q_n = F_{n+1}$ by induction.
    Again from~(b), we derive
	  \[
	    v_X(q_{n+1}-q_{n+1}(0)) = v_X(Q_n) = 1 + v_X(q_n'q_{n-1})
	    = v_X(q_n-q_n(0)) + v_X(q_{n-1}),
	  \]
	  hence $v_X(q_n-q_n(0)) = 2^{\floor{\frac n2}}$ by induction.
	  Clearly, we can write $q_n-q_n(0) = \frac z{d_n}$
	  with $d_n\in\N$ and $z\in\Z[X]$ such that
	  $d_n$ is coprime to the content $c$ of $z$.
	  On the other hand, (a) implies
	  $(-1)^nd_n = d_n\cdot(q_n(1)-q_n(0)) = z(1) \in c\Z$.
	  Therefore, $c=1$.
	  \qedhere
  \end{subproof}  
\end{proof}

\begin{table}[h]
\renewcommand{\arraystretch}{1.4}
\center
\begin{tabular}{|r|l|c|}
\hline
$n$  & $\kappa_n$ & $q_n$\\
\hline
0 & $1$ & $X$\\
1 & $\frac12$ & $1-X$ \\
2 & $\frac13$ & $X^2$ \\
3 & $\frac3{10}$ & $1-3X^2+2X^3$ \\
4 & $\frac27$ & $5X^4-4X^5$ \\
5 & $\frac{161}{572}$ &
    $1+\frac12(-35X^4+28X^5+70X^6-100X^7+35X^8)$\\
6 & $\frac{24941}{89148}$ &
$\frac1{23}(3575X^8 - 5720X^9 - 6292X^{10} + 19240X^{11}
 - 14300X^{12} + 3520X^{13})$\\
\cline{3-3}
7 & \multicolumn{2}{l|}
  {$\frac{49675943612}{177918244665}$}\\
8 & \multicolumn{2}{l|}
  {$\frac{3267335346149361824147}{11711158115225119429452}$}\\
9 & \multicolumn{2}{l|}
  {$\frac{2507700451651989905962493021537936733790431031}
         {8990773234863161759100003096510729982749072312}$}\\
10& \multicolumn{2}{l|}
  {$\frac{390583621937017677187215045781161381581437854%
          10766642680982462728116470023287868511995843}
         {140048278006628885452600904137492554179859017%
          924910241263151850844470542993943699969398879}$}\\
\hline 
\end{tabular}
\caption{$\kappa_0,\dotsc\kappa_{10}$ and $q_0,\dots q_6$}
\label{tbl-kapqn}
\end{table}

\begin{table}[h]
\center
\begin{tabular}{|r|c|c|c|}
\hline
 $n$ & $\kappa_n = \area h_n$
     & $\thet_n:=\frac{\kappa_n-\kappa_{n+1}}{\kappa_{n-1}-\kappa_n}$
     & $\kappa'_n:=\frac{\kappa_n-\thet_n\kappa_{n-1}}{1-\thet_n}$ \\
\hline
 $0$ & $1.00000000000000000000$ & & \\
 $1$ & $0.50000000000000000000$ & $0.33333333333333333333$
     & $0.25000000000000000000$ \\
 $2$ & $0.33333333333333333333$ & $0.20000000000000000000$
     & $0.29166666666666666667$ \\
 $3$ & $0.30000000000000000000$ & $0.42857142857142857143$
     & $0.27500000000000000000$ \\
 $4$ & $0.28571428571428571429$ & $0.29720279720279720280$
     & $0.27967306325515280739$ \\
 $5$ & $0.28146853146853146853$ & $0.39988492368275887552$
     & $0.27863938556045289094$ \\
 $6$ & $0.27977071835599228250$ & $0.33227684207155282888$
     & $0.27892584129464787423$ \\
 $7$ & $0.27920657437653008783$ & $0.37794408664511337726$
     & $0.27886381599022319905$ \\
 $8$ & $0.27899335949547590907$ & $0.34720353947195583609$
     & $0.27887995667262631637$ \\
 $9$ & $0.27891933053410580614$ & $0.36746593866934728811$
     & $0.27887632396190421606$ \\
$10$ & $0.27889212741232722442$ & $0.35386048463306698285$
     & $0.27887722953094531460$ \\
$11$ & $0.27888250130247112315$ & $0.36282275236161039705$
     & $0.27887701998359875188$ \\
$12$ & $0.27887900873079859727$ & $0.35683379624435139475$
     & $0.27887707102387107908$ \\
$13$ & $0.27887776246319003437$ & $0.36078940911521411257$
     & $0.27887705903280182623$ \\
$14$ & $0.27887731282303594153$ & $0.35815481970563930029$
     & $0.27887706192018075817$ \\
$15$ & $0.27887715178224762000$ & $0.35989857725744438595$
     & $0.27887706123667384356$ \\
$16$ & $0.27887709382389702266$ & $0.35873932154530262253$
     & $0.27887706140035996550$ \\
$17$ & $0.27887707303195765148$ & $0.35950754355823571103$
     & $0.27887706136147039755$ \\
$18$ & $0.27887706555709860234$ & $0.35899731864060339565$
     & $0.27887706137075986477$ \\
$19$ & $0.27887706287364424648$ & $0.35933565744435908888$
     & $0.27887706136854903338$ \\
$20$ & $0.27887706190938341130$ & $0.35911105319893171159$
     & $0.27887706136907650507$ \\
$21$ & $0.27887706156310668722$ & $0.35926004095941929301$
     & $0.27887706136895087046$ \\
$22$ & $0.27887706143870329714$ & $0.35916115956025552477$
     & $0.27887706136898082869$ \\
$23$ & $0.27887706139402243131$ & & \\
\hline 
\end{tabular}
\caption{$\thru\kappa0{23}$, $\thru\thet1{22}$ and
         $\thru{\kappa'}1{22}$ rounded to $20$ decimals}
\label{tbl-kapthetn}
\end{table}

With a Python program implementing
the formulae in Theorem~\ref{thm-qn}(b),
we are able to calculate $\thru\kappa0{23}$ and $\thru q0{24}$.
Their numerators and denominators appear to grow exponentially in length.
What fits into this paper width is listed in Table~\ref{tbl-kapqn}.
More terms are provided in the OEIS entries~\cite{Mi0,Mi2}.
During our calculations, we observe that the numerators
of $\kappa_n$ and $\kappa_{n+1}$ seem to be mostly coprime,
while the numerator of $\kappa_n$ seems to \qt{almost}
divide the denominator of $\kappa_{n+1}$.

\begin{obs}
  For $n\in\N_0$, write $\kappa_n=\frac{\mu_n}{\nu_n}$
  with $\mu_n,\nu_n\in\N$ and $\gcd(\mu_n,\nu_n)=1$.
  Then
  \begin{substate}
  \item $\mu_n$ and $\mu_{n+1}$ are coprime
    for $n\in\{0,\dotsc,4,6,\dotsc,22\}$,
    while $\gcd(\mu_5,\mu_6)=7$,
  \item $\mu_n$ divides $\nu_{n+1}$
    for $n\in\{0,1,2,4,6,8,10,11,12,15,16,17,18,19,21,22\}$,
    while \mbox{$\mu_3 = 3\cdot\gcd(\mu_3,\nu_4)$},
    $\mu_5 = 7\cdot\gcd(\mu_5,\nu_6)$,
    $\mu_7 = 19^2\cdot\gcd(\mu_7,\nu_8)$,
    $\mu_9 = 37\cdot\gcd(\mu_9,\nu_{10})$,
    $\mu_{13} = 3\cdot\gcd(\mu_{13},\nu_{14})$,
    $\mu_{14} = 7\cdot\gcd(\mu_{14},\nu_{15})$ and
    $\mu_{20} = 19^2\cdot\gcd(\mu_{20},\nu_{21})$.
  \end{substate}
\end{obs}

As for the primitive polynomials in Theorem~\ref{thm-qn}(d),
their relevant coefficients appear to be non-zero 
and no three consecutive coefficients seem to have the same sign.

\begin{obs}
  Let $n\in\enum0{24}$ and $d_n\in\N$ such that
  $\sum_ic_iX^i := d_n\cdot(q_n-q_n(0)) \in \Z[X]$ is primitive,
  as in Theorem~\ref{thm-qn}(d).
  Then
  $0\notin\{c_{i-1}c_i,c_ic_{i+1}\}\not\sseq\N$
  for $2^{\floor{\frac n2}}<i<F_{n+1}$.
\end{obs}

While each $\kappa_n$ provides an upper bound
for the stribolic constant $\kappa$,
the last Corollary from~\cite{Mi1},
stated here as a proposition,
allows us to limit $\kappa$ from below.

\begin{prp}\label{prp-kappock}
  For all $n\in\N$, we have
  $\kappa_n-1+\frac{\kappa_n}{\kappa_{n-1}} < \kappa$.
\end{prp}

Applying Proposition~\ref{prp-kappock} with $n=23$
and using the left column of Table~\ref{tbl-kapthetn}
leads to the estimate
\begin{equation}\label{eqn-kapest1}
  0.2788770612338
  < \kappa_{23}-1+\frac{\kappa_{23}}{\kappa_{22}}
  < \kappa <  \kappa_{23} < 0.2788770613941.
\end{equation}
The convergence rate of the sequence $(\kappa_n)_{n\in\N_0}$
is expressed by the quotients
\[
  \thet_n:=\frac{\kappa_n-\kappa_{n+1}}{\kappa_{n-1}-\kappa_n}
  \quad\text{for }n\in\N.
\]
Comparing their values
in the middle column of Table~\ref{tbl-kapthetn}
for of $2\leq n\leq22$
gives rise to the following hypothesis.

\begin{cnj}\label{cnj-thet}
  For all $m\in\N$, we have
  $\thet_{2m} < \thet_{2m+2} < \thet_{2m+3}<\thet_{2m+1}$.
\end{cnj}

Introducing the numbers
\[
  \kappa'_n
  := \kappa_{n-1} - \frac{\kappa_{n-1}-\kappa_n}{1-\thet_n}
   = \frac{\kappa_n-\thet_n\kappa_{n-1}}{1-\thet_n}
   = \frac{\kappa_{n-1}\kappa_{n+1}-\kappa_n^2}
          {\kappa_{n-1}-2\kappa_n+\kappa_{n+1}}
  \quad\text{for }n\in\N,
\]
given in the right column of Table~\ref{tbl-kapthetn},
Conjecture~\ref{cnj-thet} would allow for an estimate
that is considerably sharper than~\eqref{eqn-kapest1}.

\begin{prp}\label{prp-kap'}
  If Conjecture~\ref{cnj-thet} holds,
  then $\kappa'_{2m-1}<\kappa<\kappa'_{2m}$ for all $m\in\N$,
  in particular,
  $0.27887706136895087 < \kappa'_{21} < \kappa
   < \kappa'_{22} < 0.27887706136898083$.
\end{prp}
\begin{proof}
  Let $m\in\N$ and set
  $\delta_n:=\kappa_n-\kappa_{n+1}$ for $n\in\N_0$.
  Presuming Conjecture~\ref{cnj-thet}, $n>2m$ entails
  $\frac{\delta_n}{\delta_{2m-1}}
   = \thet_{2m}\dotsm\thet_n > \thet_{2m}^{n-2m+1}$.
  Therefore,
  \[
    \kappa_{2m-1} - \kappa = \sum_{n=2m-1}^\infty\delta_n
    > \delta_{2m-1}\sum_{i=0}^\infty\thet_{2m}^i
    = \frac{\delta_{2m-1}}{1-\thet_{2m}},
  \]
  that is,
  $\kappa < \kappa_{2m-1} - \frac{\delta_{2m-1}}{1-\thet_{2m}}
   = \kappa'_{2m}$.
  Likewise,
  \[
    \kappa_{2m} - \kappa = \sum_{n=2m}^\infty\delta_n
    < \delta_{2m}\sum_{i=0}^\infty\thet_{2m+1}^i
    = \frac{\delta_{2m}}{1-\thet_{2m+1}},
  \]
  hence
  $\kappa'_{2m+1}
  = \kappa_{2m} - \frac{\delta_{2m}}{1-\thet_{2m+1}}
  < \kappa$.
  The estimate $\kappa'_1<\kappa$ is now also covered
  because $\kappa'_1=\frac14<\frac{11}{40}=\kappa'_3$.
\end{proof}

Iterating the procedure of Conjecture~\ref{cnj-thet} and Proposition~\ref{prp-kap'} four times,
yields
\begin{align*}
  \kappa''''_{18} &\approx 0.2788770613689750648156,\\
  \kappa''''_{19} &\approx 0.2788770613689750647749,
\end{align*}
which we conjecture to be upper resp.\ lower bounds for $\kappa$.

\end{document}